\documentclass[11pt]{article}
\usepackage{geometry} 
\usepackage{amsmath,amsthm,amssymb}                
\usepackage{graphicx, color}
\usepackage{amssymb}
\usepackage{epstopdf}
\usepackage{pdfsync}
\usepackage{mathabx}

\def\dxy{\,dx dy}
\def\dx{\,dx}

\title{A simplified counterexample to the integral representation\\ of the relaxation of double integrals}
\author{Andrea Braides
\\ 
SISSA, via Bonomea 265, Trieste, Italy}
\date{}
\begin{document} 
\maketitle
\begin{abstract}
We show that the lower-semicontinuous envelope of a non-convex double integral may not admit a representation as a double integral. By taking an integrand with value $+\infty$ except at three points (say $-1$, $0$ and $1$) we give a simple proof and an explicit formula for the relaxation that hopefully may shed some light on this type of problems. This is a simplified version of examples by Mora-Corral and Tellini, and Kreisbeck and Zappale, who characterize the lower-semicontinuous envelope via Young measures.
\end{abstract}

Double-integral functionals defined in $L^p$ spaces of the form
\begin{equation}
F(u):=\int_{\Omega\times\Omega} f(u(x)-u(y))\dxy
\end{equation}
can be treated using the direct methods of the Calculus of Variations. 
To that end, necessary and sufficient conditions for the lower semicontinuity of $F$ with respect to weak $L^p$ topologies turn out to be the convexity and lower semicontinuity of $f$, exactly as in the case of single-integral functionals (see e.g.~\cite{BMC,P}). In the case of non-convex $f$ the parallel is lost. Indeed, in \cite{MCT} it is shown that the lower-semicontinuous envelope of $F$ cannot be represented as a double integral of the same form when the function $f$ is a simple double-well potential.
The proof in \cite{MCT} relies on the representation of the relaxed functional in terms of Young measures and on the study of the optimality conditions satisfied by such measure-valued minimizers. We now give a simple explanation of the non-representability of the relaxed functional when $f$ is a double-well  potential (or rather a ``triple-well'' potential with wells in $-1$, $0$ and $1$) with ``infinite depth''; namely,
\begin{equation}\label{f}
f(z)=\begin{cases} 0 & \hbox{if } z\in\{-1,1\}\cr
1 &   \hbox{if } z=0\cr
+\infty & \hbox{otherwise.}\end{cases}
\end{equation}
For simplicity we chose $\Omega=(0,1)$. We remark that we can extend this example to everywhere finite integrands $f$ by approximation.

We note that other examples are shown in \cite{KZ} when the integrand is of the form $f(u(x),u(y))$. In that case the functionals are not invariant by translations, so the parallel with local functionals would be with integrands depending on $(u(x),\nabla u(x))$, for which lower-semicontinuity conditions are more complex \cite{DGBDM}.
 
\smallskip
We now turn to the analysis of the counterexample.

\smallskip
{\bf Characterization of the lower-semicontinuous envelope.}  Note preliminarily that a lower bound for the lower-semicontinuous envelope $\overline F$ of $F$ with respect to the weak $L^1$-convergence is
\begin{equation}
 F_0(u):=\int_{\Omega\times\Omega} f^{**}(u(x)-u(y))\dxy,
\end{equation}
where the lower-semicontinuous convex envelope of $f$ is 
\begin{equation}
f^{**}(z)=\begin{cases} 0 & \hbox{if } z\in[-1,1]\cr
+\infty & \hbox{otherwise;}\end{cases}
\end{equation}
that is,
\begin{equation}
F_0(u)=\begin{cases} 0 & \hbox{if } {\rm ess\hbox{-}sup} \, u-{\rm ess\hbox{-}inf}\, u\le 1\cr
+\infty & \hbox{otherwise.}\end{cases}
\end{equation}
This lower bound implies that $\overline F$ is finite at most on functions $u\in L^\infty(0,1)$ such that
\begin{equation}\label{esssupinf}
{\rm ess\hbox{-}sup} \, u-{\rm ess\hbox{-}inf}\, u\le 1.
\end{equation} 

Let $u\in L^\infty(0,1)$ satisfy \eqref{esssupinf}, and let $u_j$ be a sequence weakly converging to $u$ and such that $F(u_j)<+\infty$ for all $j$. Note that for fixed $j$ the function $u_j$ can take at most two values almost everywhere and these values are at distance $1$. 
Indeed by Fubini's theorem for almost all $y\in(0,1)$ we have $u_j(x)\in\{u_j(y), u_j(y)-1, u_j(y)+1\}$ for almost every $x\in(0,1)$. Hence, there exists $z^j$ such that $u_j(x)\in\{z^j,z^j-1,z^j+1\}$  for almost every $x\in(0,1)$. 
If both values $z^j-1$ and $z^j+1$ were taken on sets of positive measure,  then 
we would have $F(u_j)=+\infty$, and a contradiction.
Hence, we can suppose that there exist $z^j$ such that $u_j(x)\in\{z^j, z^j+1\}$ almost everywhere. We can assume, up to subsequences, that $z_j\to z$, so that 
\begin{equation}\label{constraint}
z\le {\rm ess\hbox{-}inf} \,u\quad\hbox{ and }\quad {\rm ess\hbox{-}sup} \, u\le z+1,
\end{equation} and that, if we let $A^j:= \{x: u_j(x)=z+1\}$, there exists $t\in[0,1]$ such that $\lim\limits_{j\to+\infty}|A^j|= t$. Hence, we obtain 
\begin{eqnarray}\label{limFj}
\lim_{j\to+\infty}F(u_j)=\lim_{j\to+\infty}(|A^j|^2+(1-|A^j|)^2)= t^2+(t-1)^2=2t^2-2t+1.
\end{eqnarray}
Note that the minimum of $t^2+(1-t)^2$ is $1\over 2$ so that \eqref{limFj} implies that $\overline F(u)\ge {1\over 2}$ for all $u$.

Since by the convergence of $\int_{(0,1)} u_j\dx$ to $ \int_{(0,1)} u\dx$ we have
\begin{equation}\label{tuz}
t=\int_{(0,1)} u\dx-z,
\end{equation}
the limit of $F(u_j)$ can be described in terms of $\int_{(0,1)} u\dx$ and $z$ only, and is independent of the particular sequence $u_j$. 

Note conversely that if $u$ and $z$ are such that \eqref{constraint} holds, then there exist $u_j$ with $u_j\in\{z,z+1\}$ and weakly converging to $u$, so that the value $2t^2+2t+1$ is achieved on this sequence with $t$ given by \eqref{tuz}.
By optimizing in $z$ we then have a description of $\overline F(u)$ as
\begin{eqnarray}\nonumber\label{charac}
\overline F(u)&=&\min \biggl\{\Bigl(\int_{(0,1)} u\dx-z\Bigr)^2+\Bigl(\int_{(0,1)} u\dx-z-1\Bigr)^2\\ \nonumber
&&\hskip4cm :z\le {\rm ess\hbox{-}inf}\, u,\  {\rm ess\hbox{-}sup} \, u\le z+1\Bigl\}\\ \nonumber
&=&\min \biggl\{2\Bigl(\int_{(0,1)} u\dx\Bigr)^2-2(2z+1)\Bigl(\int_{(0,1)} u\dx\Bigr)+2z^2+2z+1\\
&&\hskip4cm :z\le {\rm ess\hbox{-}inf} \,u,\  {\rm ess\hbox{-}sup} \, u\le z+1\Bigl\}.
\end{eqnarray}
We can make this formula more symmetric by the change of variables $w=z+{1\over 2}$, so that
\begin{eqnarray}\nonumber\label{charac-w}
\overline F(u)&=&\min \biggl\{2\Bigl(\int_{(0,1)} u\dx\Bigr)^2-4w\Bigl(\int_{(0,1)} u\dx\Bigr)+2w^2+{1\over 2}\\
&&\hskip4cm :{\rm ess\hbox{-}sup} \, u-{1\over 2}\le w\le {\rm ess\hbox{-}inf} \,u+{1\over 2}\Bigl\}.
\end{eqnarray}
Furthermore, noting that the functionals are invariant if we add a constant to $u$, replacing $u$ by $u-\int_{(0,1)} u\dx=0$ we also have 
\begin{equation}\label{charac-w-0}
\overline F(u)=\min \biggl\{2w^2+{1\over 2} :{\rm ess\hbox{-}sup} \, u-\int_{(0,1)} u\dx-{1\over 2}\le w\le {\rm ess\hbox{-}inf} \,u-\int_{(0,1)} u\dx+{1\over 2}\Bigl\}.
\end{equation}

\smallskip
{\bf Non representability of the lower-semicontinuous envelope.} 
We now prove that there exists no $g$ such that
\begin{equation}\label{none}
\overline F(u)=\int_{\Omega\times\Omega} g(u(x)-u(y))\dxy.
\end{equation}
Note that $g$ can be assumed to be even, up to replacing $g(z)$ with ${1\over 2}(g(z)+g(-z))$.

We first describe $\overline F(u)$ more precisely in some `extreme' cases. In the first one the minimization does not involve constraint \eqref{constraint}, so that $\overline F(u)={1\over 2}$. To get this, we note that if \begin{equation}\label{mezzo} 
{\rm ess\hbox{-}sup} \, u-{\rm ess\hbox{-}inf}\, u\le {1\over 2}\end{equation}
then we can take 
$$z=\int_{(0,1)} u\dx-{1\over 2},$$
and by \eqref{mezzo} we have
$$
 {\rm ess\hbox{-}sup} \,u\le {\rm ess\hbox{-}inf}\, u+{1\over 2}\le  z+1
\quad\hbox{ and }
\quad
z\le {\rm ess\hbox{-}sup}\, u-{1\over 2}\le  {\rm ess\hbox{-}inf} \,u,
$$
and $\overline F(u)={1\over 2}$ by formula \eqref{charac}.
As a particular case of a function satisfying \eqref{mezzo} we can take $u$ a constant. In this case 
\eqref{none} would give 
\begin{equation}\label{half}
g(0)={1\over 2}\,.
\end{equation}

The other `extreme' case is when only one $z$ is involved in the minimization in \eqref{charac}; which is the case when $ {\rm ess\hbox{-}sup} \, u-{\rm ess\hbox{-}inf}\, u=1$, so that $z={\rm ess\hbox{-}inf}\,u$ and $z+1={\rm ess\hbox{-}sup} \, u$. The value of $\overline F(u)$ is then 
simply 
$$
\overline F(u)= \Bigl(\int_{(0,1)} u\dx-{\rm ess\hbox{-}inf}\, u\Bigr)^2+
\Bigl(\int_{(0,1)} u\dx-{\rm ess\hbox{-}sup}\, u\Bigr)^2\,.
$$
This can be applied, for fixed $t\in(0,1)$,  with $u$ given by
$$
u(x)= \begin{cases} 1 & \hbox{ if }x\le t\cr
0 &  \hbox{ if }x> t,\end{cases}
$$
for which $
\overline F(u)= 2t^2-2t+1$.
If \eqref{none} held true then by 
\eqref{half} we would also have
$$
\overline F(u)= (2t^2-2t+1)g(0)+ 2t(1-t) g(1)= {1\over 2}(2t^2-2t+1)+ 2t(1-t) g(1),
$$
 which would give
$$
g(1)={2t^2-2t+1\over 4t(1-t)}=  {1\over 4}\Bigl({t\over 1-t} +{1-t\over t}\Bigr).
$$
Taking different values for $t\in(0,1)$ we get different values for $g(1)$, which is a contradiction.
\smallskip

{\bf Conclusions and remarks.} Formula \eqref{charac} shows that $\overline F(u)$ is obtained by functions $u_j$ weakly converging to $u$ and oscillating between two values $z$ and $z+1$ maximizing the measure of the subset of points $(x,y)\in\Omega\times\Omega$ such that $u_j(x)=z$ and $u_j(y)=z+1$. This operation depends only on $z$, which satisfies some constraints due to the convergence of $u_j$ to $u$; minimizing the outcome in $z$ gives 
the optimal choice of $u_j$. Minimization in $z$ is unconstrained if $ {\rm ess\hbox{-}sup} \, u-{\rm ess\hbox{-}inf}\, u\le {1\over 2}$, while it is limited to a single $z$ when $ {\rm ess\hbox{-}sup} \, u-{\rm ess\hbox{-}inf}\, u=1$.
The dependence on the quantity $ {\rm ess\hbox{-}sup} \, u-{\rm ess\hbox{-}inf}\,u$ highlights the nonlocality of the recovery sequences. An example of this fact is obtained by considering constants $u=c$, for which we have minimizing sequences oscillating between $c-{1\over 2}$ and $c+{1\over 2}$, while this is not true for piecewise-constant functions: if $u$ takes only two values at distance $1$ then a recovery sequence is $u$ itself, without oscillations.

We remark that from this example we also obtain examples with finite integrand. Indeed, if $f_n$ is a sequence of functions increasingly converging to $f$ given by \eqref{f} and 
\begin{equation}
F_n(u):=\int_{\Omega\times\Omega} f_n(u(x)-u(y))\dxy,
\end{equation}
then the lower-semicontinuous envelopes $\overline F_n$ converge to $\overline F$. If there existed (convex) functions $g_n$ such that 
\begin{equation}
\overline F_n(u):=\int_{\Omega\times\Omega} g_n(u(x)-u(y))\dxy,
\end{equation}
then this would hold also for $\overline F$.

\smallskip
{\bf Acknowledgmnents.} The content of this work is a lesson of the course ``Local and nonlocal variational problems in Sobolev spaces" held at SISSA in the Winter Semester 2022-2023. The author gratefully acknowledges valuable comments by Carolin Kreisbeck.

\end{document}